\documentclass[draft]{amsart} 
\usepackage{amsmath,amssymb,amsfonts}
\usepackage{mathrsfs} 
\usepackage[dvips,final]{graphicx}

\renewcommand{\le}{\leqslant}
\renewcommand{\ge}{\geqslant}
\newtheorem{thm}{Theorem}
\newtheorem{lem}{Lema}
\newtheorem{defn}{Definition}
\newtheorem{cor}{Corrolary}
\newtheorem{rem}{Remark}

\title{Recursive state estimation for noncausal discrete-time 
descriptor systems under uncertainties}
\author{Serhiy M.Zhuk}
\address{Taras Shevchenko Kyiv National University, Faculty of Cybernetics, 
System Analysis Laboratory, Ukraine}
\email{Serhiy.Zhuk@gmail.com}
\begin{document}


\maketitle

\begin{abstract}
This paper describes a method for the online state estimation of systems 
described by a general class of linear noncausal time-varying difference descriptor equations subject to uncertainties. The method is based on the notions of a \emph{linear minimax estimation} and an \emph{index of causality} introduced here for singular difference equations. %
\emph{The online minimax estimator} is derived by the application of the dynamical programming and Moore's pseudoinverse theory to the minimax estimation problem. It coincides with Kalman's filter for regular systems. 
A numerical example of the state estimation for 2D noncasual 
descriptor system is presented.
\end{abstract}


\section{Introduction}

There is a number of physical and engineering objects most naturally modelled as systems of differential and algebraic equations (DAEs) or descriptor systems: 
microwave circuits~\cite{fav}, 
flexible-link planar parallel platforms~\cite{mills} 
and image recognition problems (noncasual image modelling)~\cite{imod}. 
DAEs arise in economics~\cite{dai}. Also nonlinear differential-algebraic 
systems are studied with help of DAEs via linearization: a batch chemical reactor model~\cite{chem}.

On the other hand there are many papers devoted 
to the mathematical processing of data obtained from the measuring device during an experiment. 
In particular, a problem of the observer design for discrete-time descriptor 
systems was studied in the \cite{dai}-\cite{twim}, the guaranteed state estimation for a linear dynamical systems was investigated in the \cite{ogn}. 
In the \cite{nikh1} the authors derive a so-called "3-block" form for the optimal filter and a corresponding 3-block Riccati equation using a maximum likelihood approach. A filter is obtained for a general class of time-varying descriptor models. The measurements are supposed to contain a noise with Gauss'es distribution. 
The obtained recursion is stated in terms of the 3-block matrix pseudoinverse. \\
In the \cite{ishixara2} the filter recursion is represented in terms of a deterministic data fitting problem solution. The authors introduce an explicit form of the 3-block matrix pseudoinverse for a descriptor system with a special structure: so their 
filter coincides with obtained in the \cite{nikh1}.

In this paper we study an observer design problem for a general class of linear noncasual time-varying descriptor models with no restrictions to a system structure. Suppose we are given an exact mathematical model of some real process and the vector $x_k$ describes the system output at the moment $k$ in the corresponding state space of the system. Also the successive measurements $y_0\dots y_k\dots$ of 
the system output $x_k$ are supposed to be available with the noise $
g_0\dots g_k\dots$ of an uncertain nature\footnote{For instance we do not have a-priory information about its distribution.}. 
Further assume that the system input $f_k$, start point $q$ and noise $g_k$ are arbitrary elements of the given set $G$. The aim of this paper is to design a minimax observer $k\mapsto\hat{x}_k$ that gives an online guaranteed estimation of the output $x_k$ on the basis of measurements $y_k$ and the structure of $G$. In~\cite{twim} minimax estimations were derived from the 2-point boundary value 
problem with the conditions at $i=0$ (start point) and $i=k$ (end point). 
Hence a recalculation of the whole 
history $\hat{x}_0\dots\hat{x}_k$ is required if the moment $k$ changes. Here we derive 
the observer $(k,y_k)\mapsto\hat{x}_k$ by 
applying dynamical programming methods 
to the minimax estimation problem 
similar to posed in the~\cite{twim}. 
We construct a map $\hat{x}$ that takes $(k,y_k)$ to $\hat{x}_k$ making it possible to assign a unique sequence of estimations $\hat{x}_0\dots\hat{x}_k\dots$ to the given sequence of observations $y_0\dots y_k\dots$ in the real time. A resulting filter recursion is stated in terms of the pseudoinverse of positive semi-defined $n\times n$- matrices.
\section{Minimax estimation problem}
Assume that $x_k\in\mathbb R^n$ is described by 
the equation\begin{equation}
  \label{eq:state}
  F_{k+1}x_{k+1}-C_k x_k = f_k, k=0,1,\dots,
\end{equation}
with the initial condition 
\begin{equation}
  \label{eq:strt}
  F_0x_0=q,
\end{equation}
and $y_k$ is given by 
\begin{equation}
  \label{eq:obs}
  y_k=H_k x_k+g_k, k=0,1,\dots,
\end{equation}
where $F_k,C_k$ are $m\times n$-matrices, $H_k$ is $p\times n$-matrix. Since we 
deal with a descriptor system we 
see that for any $k$ there is a set of vectors ${x_1^0\dots x_k^0}$ satisfying~\eqref{eq:state} while $f_i=0,q=0$. Thus the undefined inner influence caused by ${x_1^0\dots x_k^0}$ is possible to appear in the systems output. 
Also we suppose the initial condition $q$, input $\{f_k\}$ and noise $\{g_k\}$ to be unknown elements of the given set\footnote{Here and after $(\cdot,\cdot)$ denotes an inner product in an appropriate euclidean space, $\|x\|=
(x,x)^\frac 12$.}
\begin{equation}
  \label{eq:Ginft}
  G(q,\{f_k\},\{g_k\})=(S q,q)+
  \sum_0^{\infty}(S_k f_k,f_k)+(R_k g_k,g_k)\leqslant1
\end{equation}
where $S,S_k,R_k$ are some symmetric positive-defined weight matrices with the appropriate dimensions. 
The trick is to fix any $N$-partial sum of~\eqref{eq:Ginft} so that $(q,\{f_k\},\{g_k\})$ belongs to
\begin{equation}
  \label{eq:cnstr}
  \begin{split}
    &\mathscr G^N:=\{(q,\{f_k\},\{g_k\}):\\
    &(S q,q)+\sum_{k=0}^{N-1}(S_k f_k,f_k)+
    \sum_{k=0}^N(R_k g_k,g_k)\leqslant1\}
  \end{split}
\end{equation}
Then we derive the estimation $\hat{x}_N=v(N,y_N,\hat{x}_{N-1})$ considering a minimax estimation problem for $\mathscr G^N$. 
Lets denote by $\mathcal N$ a set of all $(\{x_k\},q,\{f_k\})$ such that~\eqref{eq:state} is held. 
The set $\mathscr G^N_y$ is said to be \emph{a-posteriori set}, where
\begin{equation}
  \label{Gy}
  \begin{split}  
    &\mathscr G^N_y:=\{\{x_k\}:(\{x_k\},q,\{f_k\})\in\mathcal N,\\
    &(q,\{f_k\},\{y_k-H_kx_k\})\in\mathscr G^N\}
\end{split}
\end{equation}
It follows from the definition that $\mathscr G^N_y$ consists of all possible $\{x_k\}$ causing an appearance of given $\{y_k\}$ while $(q,\{f_k\},\{g_k\})$ 
runs through $\mathscr G^N$. Thus, it's naturally to look for $x_N$ estimation \textbf{only} among the elements of $P_N(\mathscr G^N_y)$, 
where $P_N$ denotes the projection that takes $\{x_0\dots x_N\}$ to $x_N$.
\begin{defn}
A linear function $\widehat{(\ell,x_N)}$ is called a  
\emph{minimax a-posteriori estimation} if the following condition holds:
\begin{equation*}
  \begin{split}
    &\inf_{\{\tilde{x}_k\}\in\mathscr G^N_y}
    \sup_{\{x_k\}\in\mathscr G^N_y}|(\ell,x_N)-(\ell,\tilde x_N)|=\\
    &\sup_{\{x_k\}\in\mathscr G^N_y}|(\ell,x_N)-\widehat{(\ell,x_N)}|
  \end{split}
\end{equation*}
The non-negative number $$
\hat{\sigma}(\ell,N)=\sup_{\{x_k\}\in\mathscr G^N_y}|(\ell,x_N)-\widehat{(\ell,x_N)}|
$$ is called a 
\emph{minimax a-posteriori error} in the 
direction $\ell$. A map $$
N\mapsto I_N=\mathrm{dim}\{\ell\in\mathbb R^n:\hat{\sigma}(\ell,N)<+\infty\}
$$ is called \emph{an index of causality} for the pair of systems~\eqref{eq:state}-\eqref{eq:obs}.
\end{defn}
\label{s:mnmx}
Denote by $k\mapsto Q_k$ a recursive map that takes each $k\in\mathbb{N}$ to the matrix $Q_k$, where 
  \begin{equation}
    \label{eq:Qk}
    \begin{split}
      &Q_k=
      H'_kR_kH_k+F'_k[S_{k-1}-
      S_{k-1}C_{k-1}W_{k-1}^+C'_{k-1}S_{k-1}]
      F_k,\\
      &Q_0=F'_0SF_0+H'_0R_0H_0,
      W_k=Q_{k}+C'_{k}S_{k}C_{k}
    \end{split}
  \end{equation}
Let $k\mapsto r_k$ be a recursive map that takes each natural number $k$ to the vector $r_k\in\mathbb R^n$, where
\begin{equation}
    \label{eq:rk1}
    \begin{split}
      &r_k=F'_kS_{k-1}C_{k-1}W_{k-1}^+
      r_{k-1}+H'_kR_k y_k,\\
      &r_0=H'_0R_0y_0
    \end{split}
\end{equation}
and to each number $i\in \mathbb{N}$ assign the number $\alpha_i$, where
\begin{equation}
  \label{eq:alphak}
  \begin{split}
    &\alpha_i=\alpha_{i-1}+(R_iy_i,y_i)-(W_{i-1}^+
 r_{i-1},r_{i-1}),\\
 &\alpha_0=(Sg,g)+(R_0y_0,y_0)
  \end{split}
\end{equation}
The main result of this paper is formulated in the next theorem.
\begin{thm}[minimax recursive estimation]\label{t:mnmx}
  Suppose we are given a natural number $N$ and a vector $\ell\in\mathbb R^n$. Then a necessary and sufficient condition for a minimax a-posteriori error $\hat{\sigma}(\ell,N)$ to be finite is that
  \begin{equation}
    \label{eq:erfn}
 Q_N^+Q_N\ell=\ell
  \end{equation}
Under this condition we have 
\begin{equation}\label{eq:err}
    \hat{\sigma}(\ell,N)=
    [1-\alpha_N+(Q^+_Nr_N,r_N)]^\frac 12(Q^+_N\ell,\ell)^\frac 12
  \end{equation}
 and 
   \begin{equation}
     \label{eq:est}
     \widehat{(\ell,x_N)}=(\ell,Q_N^+r_N)
   \end{equation}  
\end{thm}
\begin{cor}
The index of causality $I_N$ for the pair 
of systems~\eqref{eq:state}-\eqref{eq:obs} 
can be represented as $
I_N=\mathrm{rank}(Q_N)$. 
\end{cor}
\begin{cor}[minimax obsever]
\label{c:mnmxobs}
  The online \emph{minimax observer} is given by $
  k\mapsto \hat{x}_k=Q_k^+r_k$ and
  \footnote{We assume here that $\frac 10=+\infty$.}
  \begin{equation}
    \label{eq:mobs}
    \begin{split}
      &\hat{\rho}(N)=  \min_{\{x_k\}\in\mathscr G^N_y}\max_{\{\tilde{x}_k\}\in\mathscr G^N_y}\|x_N-\tilde{x}_N\|^2=\\
      &\frac{[1-\alpha_N+(Q_N\hat x_N,\hat x_N)]}{\min_i\{\lambda_i(N)\}}
    \end{split}
  \end{equation}
where $\lambda_i(N)$ are eigen values of $Q_N$. 
In this case all possible realisations of ~\eqref{eq:state} state vector $x_N$ fill the ellipsoid $ P_N(\mathscr G^N_y)\subset\mathbb R^n$, where 
\begin{equation}
  \label{eq:PnGy}
  P_N(\mathscr G^N_y)=\{x:(Q_Nx,x)-2(Q_N\hat x_N,x)+\alpha_N\leqslant1\}
\end{equation}
\end{cor}
\begin{rem}
  If $\lambda_{min}(H'_kR_kH_k)$ grows for $k=i,i+1,\dots$ then the minimax estimation error $\hat{\rho}(k)$ becomes smaller causing $\hat{x}_k$ to get closer to thereal state vector $x_k$.
\end{rem}
In~\cite{ishixara2} Kalman's filtering problem for descriptor systems was investigated from the deterministic point of view. Authors recover Kalman's recursion to the time-variant descriptor system by a deterministic least square fitting problem over the entire trajectory: 
find a sequence $\{\hat{x}_{0|k},\dots,\hat{x}_{k|k}\}$ that minimises the following fitting error cost
\begin{equation*}
  \begin{split}
    &\mathrm{J}_k(\{x_{i|k}\}_{0}^k)=
\|F_0x_{0|k}-g\|^2+\|y_0-H_0x_{0|k}\|^2+\\
&\sum_{i=1}^k\|F_{i}x_{i|k}-C_{i-1} x_{i-1|k}\|^2+\|y_i-H_ix_{i|k}\|^2
  \end{split}
\end{equation*}
assuming that the $\mathop{\mathrm{rank}}
  \begin{smallmatrix}
    F_k\\H_k
  \end{smallmatrix}\equiv n$. 
According to~\cite[p.8]{ishixara2} the successive optimal estimates $\{\hat{x}_{0|k},\dots,\hat{x}_{k|k}\}$ resulting from the minimisation of $\mathrm{J}_k$ can be found from the recursive algorithm
  \begin{equation}
    \label{eq:fltr:r}
    \begin{split}
     &\hat{x}_{k|k}=
   P_{k|k}F'_k
   (E+C_{k-1}P_{k-1|k-1}C'_{k-1})^{-1}
   C_{k-1}\hat{x}_{k-1|k-1}\\
   &+P_{k|k}H'_k R_{k}y_k,
   \hat{x}_{0|0}=P_{0|0}(F'_0q+H'_0y_0),\\
  &P_{k|k}=\bigl(F'_k(E+C_{k-1}P_{k-1|k-1}C'_{k-1})^{-1}F_k+    H'_kH_k\bigr)^{-1},\\
  &P_{0|0}=(F'_0F_0+H'_0H_0)^{-1}
    \end{split}
  \end{equation}
\begin{cor}[Kalman's filter recursion]
\label{l:eqfltr}
  Suppose the $\mathop{\mathrm{rank}}
  \begin{smallmatrix}
    F_k\\H_k
  \end{smallmatrix}\equiv n$, and let $k\mapsto r_k$ be a recursive map that takes each 
natural number $k$ to the vector $r_k\in\mathbb R^n$, where \begin{equation}
    \label{eq:rk2}
    \begin{split}
      &r_k=H'_k y_k+F'_k C_{k-1}(C'_{k-1}C_{k-1}+Q_{k-1})^+_{k-1}r_{k-1},\\
      &r_0=F'_0q+H'_0y_0
    \end{split}
\end{equation}
Then $Q^+_kr_k=\hat{x}_{k|k}$ for each $k\in\mathbb{N}$, where $\hat{x}_{k|k}$ is given by~\eqref{eq:fltr:r} and $I_k=n$.  
\end{cor}
\section{Example}
Let us set $H_0=\left[\begin{smallmatrix}
  \frac 6{10}&&\frac{96}{100}&&0\\1000&&2\frac 3{10}&&0\\1&&\frac 1{10}&&0\\
0&&0&&0
\end{smallmatrix}\right]$,$
F_k=\left[
\begin{smallmatrix}
  1&&0&&0\\0&&1&&0
\end{smallmatrix}\right]$,$$
C_k\equiv\left[
\begin{smallmatrix}
  \frac 1{40}&&\frac 12&&0\\\frac 1{10}&&\frac 14&&\frac 3{10}
\end{smallmatrix}\right],H_k\equiv\left[ 
\begin{smallmatrix}
  k*\frac 6{10}&&k&&0\\100k&&\frac k{100}&&0\\0&&0.005&&150k*q(k)\\0.05&&10k&&0
\end{smallmatrix}\right],
$$ where $q(k)=1$ if $k$ is odd and otherwise $q(k)=0$. We derive the output $x_k$ of \eqref{eq:state} and $y_k$ assuming $f_k,g_k$ to be bounded vector-functions on the whole real axis. 
Also we set $
R_k=\mathrm{diag}\{\frac 1{11(k+1},\frac 1{22(k+1)},\frac 1{33(k+1)},\frac 1{44(k+1)}\}$, $S_k=\mathrm{diag}\{\frac 1{35(k+1)},\frac 1{70(k+1)}\}$, $S=\mathrm{diag}\{\frac 1{60},\frac 1{120}\}
$. 

\begin{figure}[h]\centering
\begin{minipage}[c]{550pt} 
\includegraphics[viewport=1 150 700 700,width=400pt,clip,angle=270]{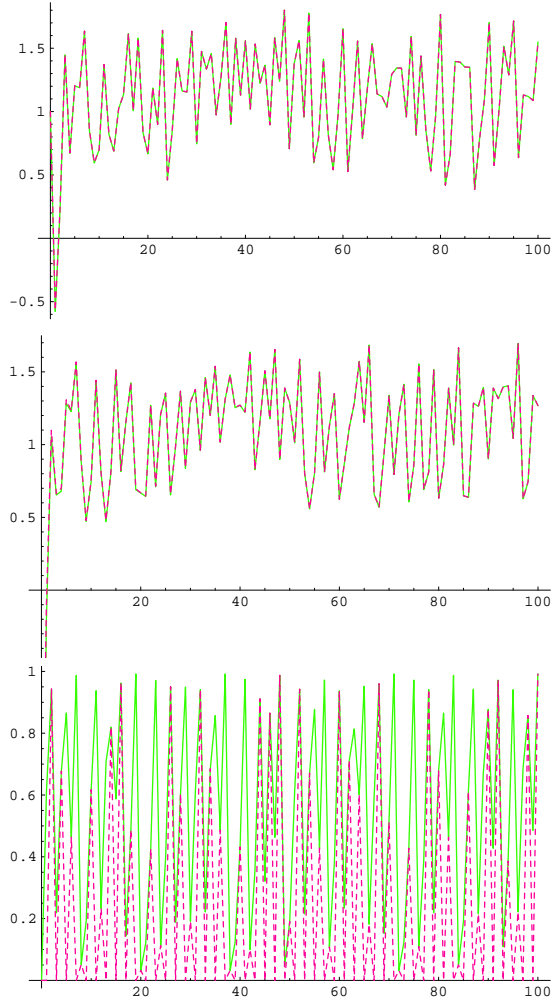}
\end{minipage}
\caption{$N=100$, output $x_{i,k}$ (solid) and observer $\hat{x}_{i,k}$ (dashed) to the left; }
\label{fig:1}
\end{figure}

\begin{figure}[h]\centering
\begin{minipage}[c]{550pt} 
\includegraphics[viewport=1 150 700 700,width=400pt,clip,angle=270]{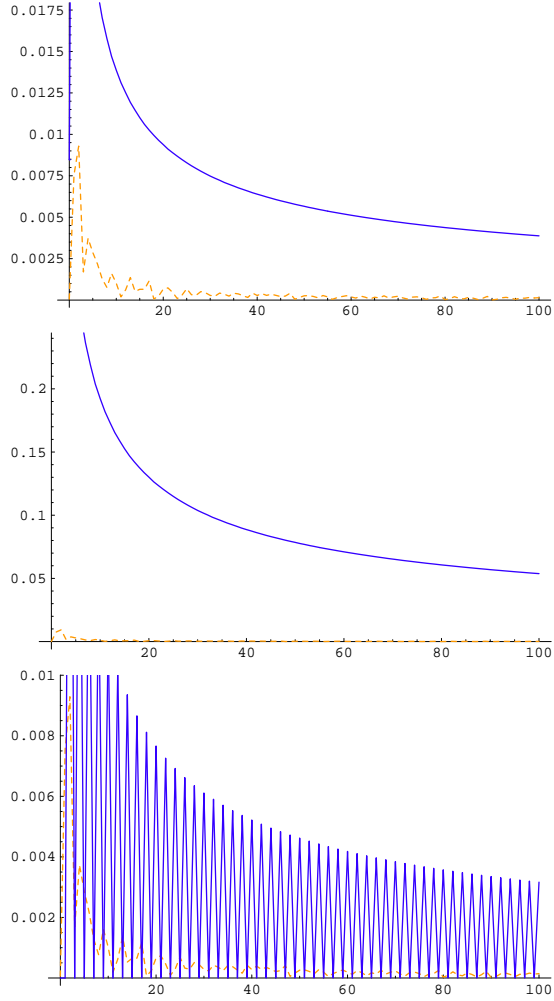}
\end{minipage}
\caption{$N=100$, real estimation error $|x_{i,k}-\hat{x}_{i,k}|$(dashed) and minimax error $\hat{\sigma}(e_i,k)$ (solid) to the right.}
\label{fig:2}
\end{figure}

We derive $\hat{x}_k$ from~\eqref{eq:rk1} and $\hat{\sigma}(e_i,k)$ from~\eqref{eq:err}, $e_i$ -- i-ort. Note that the $\mathop{\mathrm{rank}}
  \begin{smallmatrix}
    F_{2k+1}\\H_{2k+1}
  \end{smallmatrix}<3$ and $I_{2k}=3$, $I_{2k+1}<3$. Thus $\hat{x}_{3,2k+1}=0$,$$
[1-\alpha_{2k+1}+(Q^+_{2k+1}r_{2k+1},r_{2k+1})]^\frac 12(Q^+_{2k+1}\ell,\ell)^\frac 12=0
$$
but $|x_{3,2k+1}-\hat{x}_{3,2k+1}|>0$. The dynamics of $x_{i,k},\hat{x}_{i,k}$,  $|x_{i,k}-\hat{x}_{i,k}|$ and $\hat{\sigma}(e_i,k)$ is described by figures~\ref{fig:1}-\ref{fig:2}.
\par\noindent\textbf{Acknowledgements}
  It is a pleasure to thank Prof. A.Nakonechniy and Dr. V.Pichkur for insightful discussions about the key ideas presented in this paper.
\appendix
\section{Proofs.}
\begin{proof}{Proof of Theorem~\ref{t:mnmx}.}
By definition, put \begin{equation*}
  \mathbb H =\left(\begin{smallmatrix}
      H_0&&0_{pn}&&\dots&&0_{pn}\\
      0_{pn}&&H_1&&\dots&&0_{pn}\\
      \vdots&&\vdots&&\dots&&\vdots\\
      0_{pn}&&0_{pn}&&\dots&&H_N
      \end{smallmatrix}\right),
  \mathbb F =\left(\begin{smallmatrix}
      F_0&&0_{mn}&&0_{mn}&&\dots&&0_{mn}&&0_{mn}\\
      -C_0&&F_1&&0_{mn}&&\dots&&0_{mn}&&0_{mn}\\
      0_{mn}&&-C_1&&F_2&&\dots&&0_{mn}&&0_{mn}\\
      \vdots&&\vdots&&\vdots&&\dots&&\vdots&&\vdots\\
      0_{mn}&&0_{mn}&&0_{mn}&&\dots&&-C_{N-1}&&F_N\end{smallmatrix}\right)
\end{equation*}
$\mathcal{X}=\left[
\begin{smallmatrix}
  x_0\\x_1\\\vdots\\x_N
\end{smallmatrix}\right],\mathcal{Y}=\left[
\begin{smallmatrix}
  y_0\\y_1\\y_2\\\vdots\\y_N
\end{smallmatrix}\right],\\mathcal{F}=\left[
\begin{smallmatrix}
  q\\f_0\\f_1\\\vdots\\f_{N-1}
\end{smallmatrix}\right]$, $
\mathcal{G}=\left[
\begin{smallmatrix}
  g_0\\g_1\\g_2\\\vdots\\g_N
\end{smallmatrix}\right]$.

By direct calculation we obtain $
  (\ell,x_N)=(\mathcal{L},\mathcal{X}),$
$$
\mathscr G^N_y=
\{\mathcal{X}:\|\mathbb F \mathcal{X}\|_1^2+
\|\mathcal{Y}-\mathbb H \mathcal{X}\|_2^2)\leqslant 1\},
$$ where $\|\mathcal{F}\|_1^2=(Sq,q)+\sum_0^{N-1}(S_kf_k,f_k)$, $\|\cdot\|_2$ is indused by $R_k$ on the same way. This implies $$
\sup_{\{x_k\}\in\mathscr G^N_y}|(\ell,x_N-\tilde x_N)|=\sup_{\mathcal{X}\in\mathscr G^N_y}|(\mathcal{L},\mathcal{X})-(\mathcal{L},\tilde{\mathcal{X}})|
$$ Denote by $\mathrm L$  the set $R{[
  \begin{smallmatrix}
   \mathbb F '&&\mathbb H ' 
  \end{smallmatrix}]}$. 
We obviously get $$
\mathcal{L}\in\mathrm L\Leftrightarrow\sup_{\mathcal{X}\in\mathscr G^N_y}|(\mathcal{L},\mathcal{X})-(\mathcal{L},\tilde{\mathcal{X}})|<+\infty
$$ 
The application of Corollary~\ref{c:RFH} yields~\eqref{eq:erfn}. Consider a vector $\mathcal{L}\in\mathrm{L}$. Clearly $$
\inf_{\mathcal{X}\in\mathscr G^N_y}(\mathcal{L},\mathcal{X})\le(\mathcal{L},\mathcal{X})\leqslant\sup_{\mathcal{X}\in\mathscr G^N_y}(\mathcal{L},\mathcal{X}),\mathcal{X}\in\mathscr G^N_y
$$ Let $c$ denotes $\frac 12(\sup_{\mathcal{X}\in\mathscr G^N_y}(\mathcal{L},\mathcal{X})+\inf_{\mathcal{X}\in\mathscr G^N_y}(\mathcal{L},\mathcal{X}))$. 
Therefore
\begin{equation*}
  \begin{split}
    &\sup_{\mathcal{X}\in\mathscr G^N_y}|(\mathcal{L},\mathcal{X})-(\mathcal{L},\tilde{\mathcal{X}})|=\\
    &\frac 12(\mathrm s(\mathcal{L}|\mathscr G^N_y)+
(\mathrm s(-\mathcal{L}|\mathscr G^N_y))+
|c-(\mathcal{L},\tilde{\mathcal{X}})|
  \end{split}
\end{equation*}
hence 
\begin{equation}
  \label{eq:hlxt}
  \begin{split}
    &\hat{\sigma}(\ell,N)=
    \frac 12(\mathrm s(\mathcal{L}|\mathscr G^N_y)+\mathrm s(-\mathcal{L}|\mathscr G^N_y)),\\
    &\widehat{(\ell,x_N)}=
    \frac 12(\mathrm s(\mathcal{L}|\mathscr G^N_y)-\mathrm s(-\mathcal{L}|\mathscr G^N_y)),
\end{split}
\end{equation}
where $\mathrm s(\cdot|\mathscr G^N_y)$ denotes the support function of $\mathscr G^N_y$. 
Clearly, $\mathscr G^N_y$ is a convex closed 
set. Hence the equality $(\mathcal{L},\tilde{\mathcal{X}})=\widehat{(\ell,x_N)}$ is 
held for some $
\tilde{\mathcal{X}}\in\mathscr G^N_y$. 
Thus, to conclude the proof we have to calculate $\mathrm s(\mathcal{L},\mathscr G^N_y )$. Let 
  \begin{equation}
  \label{Go}
  \mathscr G^N_0=\{\mathcal{X}:\|\mathbb F \mathcal{X}\|^2+\|H\mathcal{X}\|^2\leqslant\beta_N\},
\end{equation}
where $\quad\beta_N=1-\alpha_N+(Q_N^+r_N,r_N)\ge0$.
\begin{lem}\label{l:sGy}
 \begin{equation}
    \label{eq:sGy}
    \mathrm s(\mathcal{L},\mathscr G^N_y )=
    (\ell,Q_N^+r_N)+
    \mathrm{s}(\mathcal{L}|\mathscr G^N_0)
  \end{equation}
\end{lem}
It follows from the definition of $\mathscr G^N_0$ that $
\mathrm{s}(\mathcal{L}|\mathscr G^N_0)=\mathrm{s}(-\mathcal{L}|\mathscr G^N_0)$ hence~\eqref{eq:hlxt} implies $$
\widehat{(\ell,x_N)}=(\ell,Q_N^+r_N),\hat{\sigma}(\ell)=\mathrm{s}(\mathcal{L}|\mathscr G^N_0)
$$ The application of Lemma~\ref{l:sGo} completes 
the proof.
\begin{lem}\label{l:sGo}
    \begin{equation}
    \label{eq:sGo}
    \mathrm{s}(\mathcal{L}|\mathscr G^N_0)=
    \begin{cases}
      &\sqrt{\beta_N}(Q^+_N\ell,\ell)^\frac 12,[E-Q^+_NQ_N]\ell=0,\\
      &+\infty,[E-Q^+_NQ_N]\ell\ne0
    \end{cases}
  \end{equation}
\end{lem}
\end{proof}
Let $r_k$ denote $\mathbb R^n$-valued recursive map   
\begin{equation}
    \label{eq:rk}
    \begin{split}
      &r_k=F'_k(S_{k-1}-S_{k-1}C_{k-1}P^+_{k-1}
      C'_{k-1}S_{k-1})f_{k-1}+\\
      &F'_kS_{k-1}C_{k-1}W^+_{k-1}r_{k-1}+
      H'_kR_k y_k,\\
      &r_0=F'_0Sq+H'_0R_0y_0,P_k=C'_kS_kC_k+Q_k
    \end{split}
  \end{equation}
and set
\begin{equation*}
  \begin{split}
    &\mathrm{J}(\{x_k\})=
    \|F_0x_0-g\|_S^2+\|y_0-H_0x_0\|_0^2+\\
    &\sum_{k=1}^N
    \|F_{k}x_{k}-C_{k-1} x_{k-1}-f_{k-1}\|_{k-1}^2+
    \|y_k-H_kx_k\|_k^2
  \end{split}
\end{equation*}
 where $\|g\|_S^2=(Sg,g)$, $\|f_k\|^2_k=(S_kf_k,f_k)$, $\|y_i\|_i^2=(R_iy_i,y_i)$.
\begin{lem}\label{l:prj}
  Let $x\mapsto\hat{x}_k$ be a recursive map that takes any $k\in\mathbb{N}$ to $\hat{x}_k\in\mathbb R^n$, 
where 
  \begin{equation}
    \label{eq:hx}
    \begin{split}
      &\hat{x}_k=P_k^+(C'_kS_k (F_{k+1}\hat{x}_{k+1}-f_k)+r_k),\\
      &\hat{x}_N=Q_N^+ r_N,
    \end{split}
  \end{equation}
Then $$
\min_{\{x_k\}}\mathrm{J}(\{x_k\})=\mathrm{J}(\{\hat{x}_k\})
$$
\end{lem}
\begin{proof}
By definition put $\Phi(x_0):=
\|F_0x_0-g\|_S^2+\|y_0-H_0x_0\|_0^2$
$$
\Phi_i(x_{i},x_{i+1}):=
\|F_{i+1}x_{i+1}-C_{i} x_{i}-f_{i}\|_{i}^2+\|y_{i+1}-H_{i+1}x_{i+1}\|_{i+1}^2
$$ Then we obviously get 
\begin{equation}
  \label{eq:Phi_i}
  \mathrm{J}(\{x_k\})=\Phi(x_0)+\sum_{i=0}^{N-1}\Phi_i(x_{i},x_{i+1})
\end{equation}
Let us apply a modification of Bellman's method\footnote{So-called "Kievskiy venyk" method} to the nonlinear programming task $$
\mathrm{J}(\{x_k\})\to\min_{\{x_k\}}
$$ By definition put $$
\ell_1(x_1):=\min_{x_0}\{\Phi(x_0)+\Phi_0(x_0,x_1)\}
$$ Using~\eqref{eq:Qk} and~\eqref{eq:rk} one can get $$
\Phi(x_0)=(Q_0x_0,x_0)-2(r_0,x_0)+\alpha_0\ge0,
\alpha_0:=\|g\|_S^2+\|y_0\|_0^2
$$ On the other hand it's clear that 
$$
\ell_1(x_1)=\Phi(\hat{x}_0)+\Phi_0(\hat{x}_0,x_1)=(Q_1x_1,x_1)-2(r_1,x_1)+
\alpha_1\ge0,
$$ where $\hat{x}_0=P_0^+(r_0+C'_0S_0(F_1x_1-f_0))$$$
\alpha_1:=\alpha_0+\|y_1\|_1^2+\|f_0\|_0^2-
(P_0^+(r_0-C'_0S_0f_0),r_0-C'_0S_0f_0)
$$ Considering $\ell_1(x_1)$ as an induction base 
and assuming that
\begin{equation*}
  \begin{split}
    &\ell_{i-1}(x_{i-1})=
\min_{x_{i-2}}\{\Phi_{i-2}(x_{i-2},x_{i-1})+\ell_{i-2}(x_{i-2})\}=\\
&(Q_{i-1}x_{i-1},x_{i-1})-2(r_{i-1},x_{i-1})+\alpha_{i-1}
  \end{split}
\end{equation*}
now we are going to prove that 
\begin{equation}
  \label{eq:li}
  \begin{split}
    &\ell_i(x_i)=\min_{x_{i-1}}\{\Phi_{i-1}(x_{i-1},x_{i})+\ell_{i-1}(x_{i-1})\}=\\
    &(Q_ix_i,x_i)-2(r_i,x_i)+\alpha_i  
\end{split}
\end{equation}
Note that~\cite{rkflr} for any convex function $
(x,y)\mapsto f(x,y)$ 
$$
y\mapsto\min\{f(x,y)|(x,y):P(x,y)=y\},P(a,b)=b
$$ is convex. Thus taking into account the 
definition of $\ell_1(x_1)$ one can prove by 
induction that $\ell_{i-1}$ is convex and $$
\Phi_{i-1}(x_{i-1},x_{i})+\ell_{i-1}(x_{i-1})\ge0$$ Hence\footnote{The function $x\mapsto (Ax,x)-2(x,q)+c$ is convex iff $A=A'\ge0$.} $Q_{i-1}\ge 0$, the set of global minimums $\Psi_{i-1}$ 
of the quadratic function $$
x_{i-1}\mapsto \Phi_{i-1}(x_{i-1},x_{i})+
(Q_{i-1}x_{i-1},x_{i-1})-2(r_{i-1},x_{i-1})+\alpha_{i-1}
$$ is non-empty and $\hat{x}_{i-1}\in\Psi_i$, where\footnote{The vector $\hat{x}_{i-1}$ 
has the smallest norm among other points of theminimum.} $$
\hat{x}_{i-1}=(Q_{i-1}+C'_{i-1}S_{i-1}C_{i-1})^+
(C'_{i-1}S_{i-1}(F_ix_i-f_{i-1})+r_{i-1})
$$ This implies
\begin{equation*}
  \begin{split}
    &\ell_i(x_i)=
\Phi_{i-1}(\hat x_{i-1},x_{i})+
\ell_{i-1}(\hat x_{i-1})=\\
    &
(Q_ix_i,x_i)-2(r_i,x_i)+\alpha_i ,
  \end{split}
\end{equation*}
where
\begin{equation*}
  \begin{split}
    &\alpha_i=\alpha_{i-1}+(R_iy_i,y_i)+
    (S_{i-1}f_{i-1},f_{i-1})-\\
    &(P^+_{i-1}(r_{i-1}-C'_{i-1}S_{i-1}f_{i-1}),r_{i-1}-C'_{i-1}S_{i-1}f_{i-1}),
  \end{split}
\end{equation*}
Therefore, we obtain  $$
\min_{x_N}\ell_N(x_N)=\ell_N(\hat{x}_N)=\alpha_N-(r_N,Q^+_Nr_N),
\hat{x}_N=Q^+_Nr_N
$$ so that $\min_{\{x_k\}}\mathrm{J}(\{x_k\})=
\mathrm{J}(\{\hat{x}_k\})$.
\end{proof}
\begin{cor}\label{c:RFH}
  Suppose $\mathcal{L}=[0\dots\ell]$; then $$
\mathcal{L}\in \mathscr R{[\begin{smallmatrix}
  \mathbb F '&&\mathbb H '
\end{smallmatrix}]}\Leftrightarrow
[E-Q^+_NQ_N]\ell=0
$$ and $$
\|[\begin{smallmatrix}
  \mathbb F '&&\mathbb H '
\end{smallmatrix}]^+\mathcal{L}\|^2=(Q^+_N\ell,\ell)
$$
\end{cor}
\begin{proof}
  Suppose $S_k=E,R_k=E$ for a simplicity. If $
  \mathcal{L}\in\mathscr R{[\begin{smallmatrix} \mathbb F '&&\mathbb H '\end{smallmatrix}]}$ then $$
  F'_Nz_N+H'_Nu_N=\ell,\quad F'_k z_k+H'_ku_k-C'_kz_{k+1}=0\eqno(*),
  $$ for some $z_k\in\mathbb{R}^m,u_k\in\mathbb R^p$. Let's find the projection $\{(\hat{z}_k,\hat{u}_k)\}_{k=0}^N$ of the vector $\{(z_k,u_k)\}_{k=0}^N$ onto the range of the matrix $\bigl[\begin{smallmatrix}
  \mathbb F \\\mathbb H 
\end{smallmatrix}\bigr]$. 
Lemma~\ref{l:prj} implies $$
\hat{z}_0=F_0\hat{x}_0,\hat{z}_k=F_{k}\hat{x}_k-C_{k-1}\hat{x}_{k-1},
\hat{u}_k=H_k\hat{x}_k,
\eqno(**)
$$ where 
\begin{equation*}
    \begin{split}
      &\hat{x}_k=P_k^+(C'_k F_{k+1}\hat{x}_{k+1}+r_k-C'_k z_{k+1}),
      \hat{x}_N=Q_N^+ r_N,\\
      &r_k=F'_kC_{k-1}P^+_{k-1}r_{k-1}+
      F'_k(E-C_{k-1}P^+_{k-1}C'_{k-1})z_{k}+\\
      &+H'_k u_k,
      r_0=F'_0z_0+H'_0u_0,P_k=C'_kC_k+Q_k
    \end{split}
  \end{equation*}
$(*)$ implies $r_k=C'_kz_{k+1},k=0,\dots,N-1$, $r_N=\ell$ thus $\hat{x}_N=Q^+_N\ell$, $\hat{x}_k=P_k^+C'_k F_{k+1}\hat{x}_{k+1}$ or $\hat{x}_k=\Phi(k,N)Q^+_N\ell,$$$
\Phi(k,N)=P_k^+C'_k F_{k+1}\Phi(k+1,N),\Phi(s,s)=E
$$ Combining this with $(**)$ we obtain
\begin{equation}
  \label{eq:hzhu}
  \begin{split}
    &\hat{z}_k=
(F_{k}\Phi(k,N)-C_{k-1}\Phi(k-1,N))Q^+_N\ell,\\&\hat{u}_k=
H_k\Phi(k,N)Q^+_N\ell,\hat{z}_0=F_0\Phi(0,N)Q^+_N\ell
  \end{split}
\end{equation}
By definition, put $U(0)=Q_0$,
\begin{equation*}
  \begin{split}
    &U(k)=\Phi'(k-1,k)U(k-1)\Phi(k-1,k)+\\
    &H'_kH_k+
    F_k(E-C_{k-1}P^+_{k-1}C'_{k-1})^2F_k
  \end{split}
\end{equation*}
It now follows that $$
\|[\begin{smallmatrix}
  \mathbb F '&&\mathbb H '
\end{smallmatrix}]^+\mathcal{L}\|^2=
\sum_0^N\|\hat{z}_N\|^2+\|\hat{u}_N\|^2=
(U(N)Q^+_N\ell,Q^+_N\ell)
$$ It's easy to prove by induction that $Q_k=U(k)$.

Since $$
\mathcal{L}\in R{[\begin{smallmatrix}
  \mathbb F '&&\mathbb H '
\end{smallmatrix}]}
$$ we obtain by substituting $\hat z_k,\hat u_k$ into $(*)$ $$
F'_N\hat{z}_N+H'_N\hat{u}_N=\ell
$$ On the other hand~\eqref{eq:Qk} and~\eqref{eq:hzhu} imply $$
 F'_N\hat{z}_N+H'_N\hat{u}_N=\ell\Rightarrow
 [E-Q^+_NQ_N]\ell=0
$$ 
Suppose that $[E-Q^+_NQ_N]\ell=0$. 
To conclude the proof we have to show that 
$$
(\ell,x_N)=(Q^+_N\ell,Q_Nx_N)=0, \forall[x_0\dots x_N]\in
\mathscr N {[\begin{smallmatrix}\mathbb F ,\mathbb H \end{smallmatrix}]}
$$ By induction, fix $N=0$. If $F_0x_0=0,H_0x_0=0$, then $Q_0x_0=0$. We say that $[x_0\dots x_{k}]\in\mathscr N {[\begin{smallmatrix}\mathbb F ,\mathbb H \end{smallmatrix}]}$ if $$
F_0x_0=0,H_0x_0=0,F_sx_s=C_{s-1}x_{s-1},H_sx_s=0,
$$
Suppose $Q_{k-1}x_{k-1}=0,\forall[x_0\dots x_{k-1}]\in\mathscr N {[\begin{smallmatrix}\mathbb F ,\mathbb H \end{smallmatrix}]}$ and fix any $[x_0\dots x_k]\in
\mathscr N {[\begin{smallmatrix}\mathbb F ,\mathbb H \end{smallmatrix}]}
$. Then $F_kx_k=C_{k-1}x_{k_1},H_kx_k=0$. Combining this with~\eqref{eq:Qk} we obtain $$
Q_kx_k=F'_k(E-C_{k-1}P^+_{k-1}C'_{k-1})C_{k-1}x_{k-1}\eqno(*)
$$ We show that $Q_k\ge0$ in the proof of Theorem~\ref{t:mnmx}. 
One can see that 
\begin{equation*}
  \begin{split}
    &\left[\begin{smallmatrix}
  C_{k-1}\\Q^{\frac 12}_{k-1}
\end{smallmatrix}\right]^+=\\
&[(C'_{k-1}C_{k-1}+Q_{k-1})^+C'_{k-1},(C'_{k-1}C_{k-1}+Q_{k-1})^+Q^{\frac 12}_{k-1}]
  \end{split}
\end{equation*}
Since $$
[\begin{smallmatrix}
  C_{k-1}\\Q^{\frac 12}_{k-1}
\end{smallmatrix}][\begin{smallmatrix}
  C_{k-1}\\Q^{\frac 12}_{k-1}
\end{smallmatrix}]^+[\begin{smallmatrix}
  C_{k-1}\\Q^{\frac 12}_{k-1}
\end{smallmatrix}]x_{k-1}=[\begin{smallmatrix}
  C_{k-1}\\Q^{\frac 12}_{k-1}
\end{smallmatrix}]x_{k-1}
$$ we obviously get $$
C_{k-1}P^+_{k-1}C'_{k-1}C_{k-1}x_{k-1}=C_{k-1}x_{k-1}\Rightarrow Q_kx_k=0
$$ as it follows from $(*)$. 
This completes the proof.
\end{proof}
\begin{proof}{Proof of Lemma~\ref{l:sGy}.}
Taking into account the definitions of the 
matrices $\mathbb F ,\mathbb H $ and \eqref{Gy} we clearly have $$
\mathscr G^N_y=\{\mathcal{X}:\|\mathbb F \mathcal{X}\|^2+\|\mathcal{Y}-\mathbb H \mathcal{X}\|^2\leqslant 1\}
$$ Let $\hat{\mathcal{X}}$ be a minimum of the 
quadratic function $
\mathcal{X}\mapsto\|\mathbb F \mathcal{X}\|^2+\|\mathcal{Y}-\mathbb H \mathcal{X}\|^2$. 
It now follows that $$
\mathscr G^N_y=\hat{\mathcal{X}}+\mathscr G^N_0
\Rightarrow \mathrm s(\mathcal{L}|\mathscr G^N_0)=(\mathcal{L},\hat{\mathcal{X}})+\mathrm{s}(\mathcal{L}|\mathscr G^N_0)
$$ The application of Lemma~\ref{l:prj} yields $$
(\mathcal{L},\hat{\mathcal{X}})=(\ell,Q_N^+ r_N)
$$ This completes the proof.
\end{proof}
\begin{proof}{Proof of Lemma~\ref{l:sGo}.}
Suppose the function $f:\mathbb R^n\to R^1$ is convex and closed. Then~\cite{rkflr} the support function $s(\cdot|\{x:f(x)\leqslant0\})$ of the set $\{x:f(x)\leqslant0\}$ 
is given by
$$ 
s(z|\{x:f(x)\leqslant0\})=\mathrm{cl}\inf_{\lambda\ge0}\{\lambda f^*(\frac z\lambda)\}
$$ To conclude the proof it remains to compute the support function of $\mathscr G^N_0$ according to this rule and then apply Corollary~\ref{c:RFH}. 
\end{proof}
\begin{proof}{Proof of Corollary~\ref{l:eqfltr}.}
The proof is by induction on $k$. For $k=0$, there is nothing to prove. The induction hypothesis is $
P_{k-1|k-1}=Q_{k-1}^{-1}$. 
  Suppose $S$ is $n\times n$-matrix such that $S=S'>0$, $A$ is $m\times n$-matrix; then 
  \begin{equation}
    \label{eq:AS}
    A(S^{-1}+A'A)^{-1}=(E+ASA')^{-1}AS
  \end{equation}
Using~\eqref{eq:AS} we get 
\begin{equation}
  \label{eq:AS1}
ASA'=[E+ASA']A[A'A+S^{-1}]^{-1}A'  
\end{equation}

Combining \eqref{eq:AS1} with the induction 
assumption we get the following 
\begin{equation*}
  \begin{split}
    &E+C_{k-1}P_{k-1|k-1}C'_{k-1}=\\
    &E+[E+C_{k-1}P_{k-1|k-1}C'_{k-1}]\times\\
&\times C_{k-1}
[Q_{k-1}+C'_{k-1}C_{k-1}]^{-1}C'_{k-1}
  \end{split}
\end{equation*}
By simple calculation from the previous 
equality follows 
\begin{equation*}
  \begin{split}
    &E-C_{k-1}(Q_{k-1}+C'_{k-1}C_{k-1})^{-1}C'_{k-1}=
    \\
    &(E+C_{k-1}P_{k-1|k-1}C'_{k-1})^{-1}
  \end{split}
\end{equation*}
Using this and~\eqref{eq:Qk},\eqref{eq:fltr:r} 
we obviously get $Q_k^{-1}=P_{k|k}$. 

It follows from the definitions that $Q^{-1}_0r_0=\hat{x}_{0|0}$. Suppose that $
Q^{-1}_{k-1}r_{k-1}=\hat{x}_{k-1|k-1}$. 
The induction hypothesis and~\eqref{eq:AS} 
imply
\begin{equation*}
  \begin{split}
    &(E+C_{k-1}P_{k-1|k-1}C'_{k-1})^{-1}C_{k-1}\hat{x}_{k-1|k-1}=\\
    &C_{k-1}(C'_{k-1}C_{k-1}+Q_{k-1})^{-1}_{k-1}r_{k-1}
  \end{split}
\end{equation*}
Combining this with \eqref{eq:fltr:r}, \eqref{eq:rk2} and using $Q_k^{-1}=P_{k|k}$ we obtain $$
\hat{x}_{k|k}=Q_k^{-1}(F'_kC_{k-1}(C'_{k-1}C_{k-1}+Q_{k-1})^+_{k-1}r_{k-1}+H'_ky_k)
$$ This concludes the proof.
\end{proof}
\begin{proof}{Proof of Corollary~\ref{c:mnmxobs}.}
If $I_k<n$ then $\mathrm{rank}(Q_)<n$ hence $\lambda_{min}(Q_k)=0$. In this case there is a direction $\ell\in\mathbb R^n$ such that $\hat\sigma(\ell,k)=+\infty$. So $\hat\rho(k)=+\infty$. 
If $I_k=n$ then we clearly have
 \begin{equation*}
   \begin{split}
     & \min_{\{x_k\}\in\mathscr G^N_y}\max_{\{\tilde{x}_k\}\in\mathscr G^N_y}\|x_N-\tilde{x}_N\|^2=\\
      &\min_{\{x_k\}\in\mathscr G^N_y}\max_{\{\tilde{x}_k\}\in\mathscr G^N_y}
      \{\max_{\|\ell\|=1}|(\ell,x_N-\tilde{x}_N)|\}^2=\\
     &\{\min_{\mathscr G^N_y}\max_{\|\ell\|=1}
     \max_{\{\tilde{x}_k\}\in\mathscr G^N_y}|(\ell,x_N-\tilde{x}_N)|\}^2\ge\\
     &\{\max_{\|\ell\|=1}\min_{\{x_k\}\in\mathscr G^N_y}
     \max_{\{\tilde{x}_k\}\in\mathscr G^N_y}|(\ell,x_N-\tilde{x}_N)|\}^2=\\
     &[1-\alpha_N+(Q^+_Nr_N,r_N)]\max_{\|\ell\|=1}(Q^+_N\ell,\ell)=\\
     &\frac{[1-\alpha_N+(Q^+_Nr_N,r_N)]}{\min_i\{\lambda_i(N)\}}
   \end{split}
 \end{equation*}
On the other hand Theorem~\ref{t:mnmx} implies
\begin{equation*}
  \begin{split}
    &\max_{\{\tilde{x}_k\}\in\mathscr G^N_y}\|\hat{x}_N-\tilde{x}_N\|^2=\\
&\{\max_{\|\ell\|=1}\max_{\{\tilde{x}_k\}\in\mathscr G^N_y}|(\ell,x_N-\tilde{x}_N)|\}^2=\\
&\{\max_{\|\ell\|=1}[1-\alpha_N+(Q^+_Nr_N,r_N)]^\frac 12(Q^+_N\ell,\ell)^\frac 12\}^2
  \end{split}
\end{equation*}
It follows now from $I_N=n$ that $\mathscr G^N_y$ is a 
bounded set.

The equality $I_N=n$ implies $[E-Q^+_NQ_N]=0$ for 
a given $N$. 
It follows from Lemmas~\ref{l:sGy},\ref{l:sGo} 
that
\begin{equation}\label{eq:zrch}
  \begin{split}
    &s(\ell|P_N(\mathscr G^N_y))=s(P'_N\ell|\mathscr G^N_y)=s(\mathcal{L}|\mathscr G^N_y)=\\
    &(\ell,Q_N^+r_N)+\sqrt{\beta_N}(Q^+_N\ell,\ell)^\frac 12
  \end{split}
\end{equation}
for any $\ell\in\mathbb R^n$. By 
Young's theorem~\cite{rkflr}, \eqref{eq:zrch}, so that
\begin{equation*}
  \begin{split}
    &P_N(\mathscr G^N_y)=
\{x\in\mathbb R^n:(x,\ell)\leqslant
s(\ell|P_N(\mathscr G^N_y)),
\forall\ell\in\mathbb R^n\}=\\
    &\{x\in\mathbb R^n:\sup_\ell\{(x,\ell)-(\ell,\hat{x}_N)-
    \sqrt{\beta_N}(Q^+_N\ell,\ell)^\frac 12\}\leqslant0\}=\\
    &\{x\in\mathbb R^n:(Q_Nx,x)-2(Q_N\hat x_N,x)+\alpha_N\leqslant1\}
  \end{split}
\end{equation*}
\end{proof}

\end{document}